\documentclass[12pt]{amsart}
\usepackage[russian]{babel}
\usepackage{amsmath,amsfonts,amssymb,euscript,eufrak}

\theoremstyle{plain}

\newtheorem{Proposition}{Утверждение}
\newtheorem{Theorem}{Теорема}

\newtheorem{Note}{Замечание}

\def\Dom{{\frak D}}

\def\al{\alpha}
\def\la{\lambda}

\def\dist{\operatorname{dist}}

\def\eps{\varepsilon}

\def\\singsupp{\operatorname{\sing\ supp}}
\def\sign{\operatorname{sign}}

\renewcommand\Re{\mathrm{Re}\,}
\renewcommand\Im{\mathrm{Im\,}}

\def\bR{{\Bbb R}}

\def\bC{{\Bbb C}}
\def\bN{{\Bbb N}}
\def\L{{\EuScript L}}
\def\R{{\mathcal R}}
\def\t{{\mathfrak t}}

\def\h{{\mathfrak h}}
\def\w{{\mathrm w}}

\def\({\left(}
\def\){\right)}

\def\U{\mathrm{U}}
\def\Ai{\mathrm{Ai}}
\def\Bi{\mathrm{Bi}}

\title{Спектральные свойства комплексного оператора Эйри на полуоси}

\bigskip

\author{Савчук А.~М., Шкаликов А.~А.}
\thanks{Работа выполнена при поддержке РФФИ, грант No 16-01-00706.}

\begin{document}

\maketitle

\begin{abstract}
В работе доказана теорема о полноте системы корневых функций оператора Шредингера $L=-d^2/dx^2+p(x)$ на полуоси $\bR_+$
с потенциалом $p$, при котором оператор $L$ оказывается максимально секториальным. Применение этой теоремы к оператору
Эйри $\L_c=-d^2/dx^2+cx$, $c=const$, влечет за собой полноту системы собственных функций этого оператора в случае
$|\arg c|<2\pi/3$. С использованием более тонких методов в работе доказана теорема о сохранении полноты системы
собственных функций этого специального оператора при выполнении условия $|\arg c|<5\pi/6$.
\end{abstract}

\section*{Введение}

Основное содержание этой статьи связано с изучением оператора

\begin{equation}\label{eq:main}
\L_c=-\frac{d^2}{dx^2}+cx
\end{equation}
на полуоси $x\in[0,+\infty)$ с краевым условием Дирихле в нуле. Нас будет интересовать случай невещественной константы
$c$. Основной результат будет получен в теореме 1: {\sl собственные функции этого оператора отвечают простым
собственным значениям и образуют полную систему в пространстве $L_2(\mathbb R_+)$ при условии  $|\arg c| < 5\pi/6$.}

 Мы рассмотрим также оператор
\begin{equation}\label{eq:main1}
\L_{c, \alpha}=-\frac{d^2}{dx^2} +cx^\alpha, \qquad x\in[0,\infty),
\end{equation}
где $\alpha >0,\ \,  c\in \mathbb C\setminus \mathbb R_-$. Изучение этого оператора удобно связать с операторами
более общего вида
\begin{equation}\label{eq:main2}
Ly =-y'' + p(x)y,\quad p(x) = q(x)\pm ir(x), \quad x\in [0, \infty),
\end{equation}
 где
 \begin{equation}\label{eq:cond}
r(x) \ge M_0, \ \ \ q(x) \ge c_0 r(x) +M_1, \quad \underline{\lim}_{x\to\infty} x^{-\alpha}r(x) \ge a >0, \ \ \alpha >0,
\end{equation}
 а $M_0, \, M_1,\, c_0$ ---  вещественные константы, возможно отрицательные. Функции  $q$  и $r$
достаточно считать локально суммируемыми.  Мы получим теорему 2: {\sl  если выполнены условия \eqref{eq:cond} и
оператор $L_D$ порожден дифференциальным выражением \eqref{eq:main2}  и краевым условием Дирихле в нуле, то  система
его корневых функций  полна в пространстве  $L_2(\mathbb R_+)$ при условии  $|\gamma| < 2\alpha\pi/(2+\alpha)$,  где
$\gamma = \arg(\pm i+c_0)\in (0,\pi)$. Более того, эта система образует базис для метода суммирования Абеля--Лидского.}

Заметим, что потенциал $cx^\alpha$ представим в виде $|c|(\cos\gamma +i \sin\gamma)x^\alpha$, поэтому теорема о полноте
для оператора \eqref{eq:main1} справедлива при $|\arg c| < 2\alpha\pi/(2+\alpha)$. В частности, при $c=i$  теорема о
полноте  справедлива при $\alpha > 2/3$. Теорема 2 представляется более общей, но в применении к оператору
\eqref{eq:main} из нее следует утверждение о полноте только при $|\arg c| < 2\pi/3$.  Доказательство теоремы 2
получается из общей теории, ведущей начало от работы Келдыша \cite{Ke1}. Эта теория развивалась многими авторами (см.
подробности в \cite[$\S 4$]{Shka1}). Приведенная здесь теорема есть обобщение теоремы Лидского \cite{Li1}, полученной
им  при  условии, когда в \eqref{eq:cond} постоянная $c_0=0$. Но утверждение теоремы 1 при $|\arg c| \in[2\pi/3,
5\pi/6)$ (которое не следует из теоремы 2) является существенно более тонким результатом. Он составляет наиболее важную
часть работы.

Сформулированные результаты были получена авторами в 1999 г., вскоре после обсуждения этих задач с Дависом (см. работу
\cite{Da1}). Однако авторы откладывали публикацию, надеясь решить задачу о полноте для оператора $\L_c$ полностью.
Недавно Б.~С.~Митягин обратил наше внимание  на поставленную Я.~Алмогом (Y.~Almog) проблему \cite{Alm1}: {\sl будет ли
система собственных функций оператора $\L_{i,\alpha} =-d^2/dx^2 + ix^\alpha$  полной при $\alpha \in(0, 2/3]$? } Ответа
на этот вопрос мы не знаем, но очевидной является связь этой задачи с результатом теоремы 1. В частности, мы не
сомневаемся в справедливости следующей гипотезы: {\sl найдется число $\alpha_0< 2/3$, такое, что собственные функции
$\L_{i,\alpha}$ образуют полную систему в $\L_2(\mathbb R_+)$ при $\alpha \in (\alpha_0, 2/3]$. }

Здесь мы рассматриваем операторы на полуоси $\mathbb R_+$,  но результаты сохраняются для всей оси, если потенциалы
продолжены на всю ось четным образом (см. замечание в конце работы). Обычно оператором Эйри называют оператор $\L_c$
при $c=i$ (см., например, \cite{Alm2}), но мы сохраняем это название для произвольного $c\in\mathbb C\setminus \{0\}$.
Оператор Эйри изучался в связи  с известной в гидромеханике задачей Орра--Зоммерфельда (см. работы \cite{Shka2}--
\cite{Shka3}  и имеющиеся там ссылки). Этот оператор связан также и с другими задачами механики (см., например, работу
Алмога \cite{Alm2}). Для несамосопряженных операторов важно знать не только локализацию спектра, но и иметь информацию
об $\eps$-псевдоспектре $\sigma_\eps:=\{\la\in\bC:\|(L-\la I)^{-1}\|\ge\eps^{-1}\}$. В этом направлении отметим,
например,  работы Виолы, Ембри, Зигла, Крейчирика, Трефезена, Тэйтора, Хенри \cite{TE}-- \cite{HK}.

В теме, посвященной изучению оператора Шредингера с комплексным потенциалом $p= q+ir$  можно выделить еще три
направления.  К первому направлению отнесем работы, в которых функция $r$  в некотором смысле подчинена функции $q$, и
соответствующий оператор является  возмущением самосопряженного. Здесь отметим  работы Аддучи, Виолы, Джакова, Зигла,
Митягина и Шкаликова  \cite{AM1}---\cite{MS2}. Ко второму направлению отнесем работы, в которых потенциал чисто мнимый,
т.е. $q(x)\equiv 0$.  Здесь выделим работы Дэвиса, Зигла, Крейчирика, Куидлаарса, Туманова и Шкаликова \cite{Da1},
\cite{Da2}-- \cite {He3}. К третьему направлению можно отнести работы по оператору Шредингера с так называемом
$PT$-симметричным потенциалом $p(x) =- \overline{p(-x)}$.  В этом направлении отметим работы Бендера, Ботчера, Зигла,
Крейчирика, Еременко, Габриэлова и Шапиро   \cite{BB}--\cite{EGS}. В контексте нашей работы  важно отметить работу
Гребенкова, Хеффера и Хенри \cite{GHH}, в которой оператор Эйри изучался при $c=i$ на полуоси и всей  оси. В частности,
было доказано, что собственные функции такого оператора на полуоси образуют полную систему, но не базис, а на всей оси
спектр этого оператора пуст.

Мы ограничиваемся здесь рассмотрением только краевого условия Дирихле для простоты и конкретности изложения.
Утверждения теорем и их доказательства сохраняются, если вместо условия Дирихле рассматривать условие $y'(0) +hy(0) =0,
 h\in \mathbb R$.

\section{Определение оператора $\L_c$  и его основные свойства}

Дадим более точное определение  оператора $\L_c$ с краевым условием Дирихле. А именно, считаем, что $\L_c$
 определен в  пространстве $L_2(\bR_+)$ дифференциальным выражением
\begin{equation*}
l(y)=-y''+cxy,\qquad x\in[0,+\infty),
\end{equation*}
на области
$$
\Dom(\L_c)=\{y\in L_2(\bR_+):y\in W_{2,loc}^2,\ l(y)\in L_2(\bR_+),\ y(0)=0\}.
$$
Константу $c$ будем считать комплексной, причем $c\in\bC\setminus(-\infty,0]$. Иными словами, нас будет интересовать
случай $\gamma:=\arg c\in(-\pi,\pi)$. Очевидно, оператор $\L_c$ плотно определен так как его область содержит
бесконечно дифференцируемые функции с компактным носителем на интервале $(0,+\infty)$, множество которых плотно в
$L_2(\bR_+)$.

Далее мы будем иметь дело со специальными функциями --- решениями уравнения Эйри $y''=zy$. Хорошо известно (см.,
например, \cite[Глава IV, \S 1]{Fed} и \cite[\S 10.4]{AbSt}), что это уравнение имеет пару линейно независимых решений
$\Ai(z)$ и $\Bi(z)$, с начальными условиями
\begin{gather*}
\Ai(0)=\frac1{3^{2/3}\Gamma(2/3)},\quad \Ai'(0)=-\frac1{3^{1/3}\Gamma(1/3)},\qquad
\Bi(0)=\frac1{3^{1/6}\Gamma(2/3)},\quad \Bi'(0)=\frac{3^{1/6}}{\Gamma(1/3)},\notag
\end{gather*}
причем для Вронскиана этих функций справедливо равенство
$$
 W(\Ai,\Bi)=\Ai(z)\Bi'(z)-\Ai'(z)\Bi(z)=1/\pi.\label{eq:Aiinit}
$$
Обе этих функции являются целыми функциями порядка $3/2$ и типа $2/3$. В области $|\arg z|<\pi-\eps$, где $\eps>0$
произвольно мало, при $|z|\to\infty$ функция $\Ai(z)$ допускает асимптотические представления\footnote{Здесь и далее мы
будем фиксировать ветви функций $z^\al$ условием $\arg z\in[-\pi,\pi)$.}
\begin{equation}\label{eq:Ai1}
\Ai(z)=\frac{1}{2\sqrt\pi} z^{-1/4}e^{-\frac23z^{3/2}}\left(1+O(z^{-3/2})\right),\qquad
\Ai'(z)=-\frac1{2\sqrt\pi}z^{1/4}e^{-\frac23z^{3/2}}\left(1+O(z^{-3/2})\right).
\end{equation}
Представления \eqref{eq:Ai1} можно дифференцировать произвольное количество раз. Видно, что функции $\Ai(z)$ и
$\Ai'(z)$ экспоненциально убывают при $|z|\to\infty$ на любом луче в секторе $|\arg z|<\pi/3$. Наконец, на луче $\arg
z= \pi$ для функции $\Ai(-x)$, $x>0$, справедливо представление
\begin{equation*}
\Ai(-x)=\frac1{\sqrt{\pi}}|x|^{-1/4}\left(\sin\left(\frac23|x|^{3/2}+\frac\pi4\right)+O(|x|^{-3/2})\right),\quad|x|\to\infty.
\end{equation*}
Все нули $z_k$ функции $\Ai(z)$ просты, лежат на луче $\arg z=-\pi$ и
\begin{equation}\label{eq:Aizero}
z_k=-\left[\frac32\pi\left(k-\frac14\right)\right]^{2/3}+O(k^{-4/3}),\qquad k=1,\,2,\,\dots.
\end{equation}
Вместо функции $\Bi(z)$ нам удобно будет использовать функцию $\U(z):=\Bi(z)-\sqrt3\Ai(z)$, подчиненную начальным
условиям
\begin{equation*}
\U(0)=0,\quad \U'(0)=\frac{2\cdot3^{1/6}}{\Gamma(1/3)},\qquad \text{причем}\ \ W(\Ai,\U)=1/\pi.
\end{equation*}
В секторе $|\arg z|<\pi/3$ эта функция имеет асимптотическое представление
\begin{equation}\label{eq:U1}
\U(z)=\frac1{\sqrt\pi}z^{-1/4}e^{\frac23z^{3/2}}\left(1+O(z^{-3/2})\right),\qquad
\U'(z)=\frac1{\sqrt\pi}z^{1/4}e^{\frac23z^{3/2}}\left(1+O(z^{-3/2})\right).
\end{equation}
Таким образом, функции $\U(z)$ и $\U'(z)$ экспоненциально растут при $|z|\to\infty$ на любом луче этого сектора.
Асимптотические представления функций $\U(z)$ и $\U'(z)$ в других секторах комплексной плоскости также хорошо известны,
но здесь они нам не потребуются.
\begin{Proposition}\label{pr:inverse}
Область $\Dom(\L_c)$ совпадает с множеством функций вида
\begin{equation}\label{eq:L-1}
y(x)=\pi c^{-1/3}\left[\Ai(c^{1/3}x)\int_0^x \U(c^{1/3}t)f(t)\,dt+\U(c^{1/3}x)\int_x^\infty
\Ai(c^{1/3}t)f(t)\,dt\right],
\end{equation}
где функция $f(x)$ пробегает все пространство $L_2(\bR_+)$. Для любой функции $y\in\Dom(\L_c)$:
\begin{equation}\label{eq:decr}
x^{1/2}y(x)\to 0,\quad y'(x)\to0\quad\text{при }x\to+\infty,\qquad x^{1/2}y(x),\ y'(x)\in L_2(\bR_+).
\end{equation}
Равенство \eqref{eq:L-1} задает в пространстве $L_2(\bR_+)$ ограниченный оператор, являющийся обратным к оператору
$\L_c$.
\end{Proposition}
\begin{proof}
Заметим, что луч $c^{1/3}x$, $x>0$, лежит в секторе $|\arg z|<\pi/3$ комплексной плоскости. Это означает, что функция
$\Ai(c^{1/3}t)$ экспоненциально убывает, а значит несобственный интеграл в \eqref{eq:L-1} сходится. Поскольку $f\in
L_2(\bR_+)$, то функция $y$, определенная равенством \eqref{eq:L-1}, лежит в пространстве $W_{2}^1[0,b]$ для любого
конечного $b$. Дифференцируя, получим
\begin{equation}\label{eq:L-1'}
y'(x)=\pi\left[\Ai'(c^{1/3}x)\int_0^x \U(c^{1/3}t)f(t)\,dt+\U'(c^{1/3}x)\int_x^\infty \Ai(c^{1/3}t)f(t)\,dt\right],
\end{equation}
откуда следует, что $y'\in W_2^1[0,b]$ для любого конечного $b$. Дифференцируя еще раз, получаем
\begin{gather*}
y''=\pi\left[\Ai'(c^{1/3}x)\U(c^{1/3}x)-\U'(c^{1/3}x)\Ai(c^{1/3}x)\right]f(x)+\\+\pi
c^{1/3}\left[\Ai''(c^{1/3}x)\int_0^x \U(c^{1/3}t)f(t)\,dt+\U''(c^{1/3}x)\int_x^\infty
\Ai(c^{1/3}t)f(t)\,dt\right]=-f(x)+c x y(x),
\end{gather*}
т.е.  $y\in W^2_2[0,b]$ при любом $b>0$ и $l(y)=f \in L_2(\mathbb R_+)$. Так как $y(0)=0$ то $y\in \Dom(\L_c).$
 Верно и обратное утверждение. А именно, пусть $y\in \Dom(\L_c)$ и $ l(y)= f\in L_2(\mathbb R_+)$.
  Согласно классической теореме об общем виде  решения дифференциального уравнения,
\begin{multline*}
y=C_1\Ai(c^{1/3}x)+C_2\U(c^{1/3}x)+\\+\pi c^{-1/3}\left[\Ai(c^{1/3}x)\int_0^x
\U(c^{1/3}t)f(t)\,dt+\U(c^{1/3}x)\int_x^\infty \Ai(c^{1/3}t)f(t)\,dt\right],
\end{multline*}
где  $C_1$, $C_2$ --- постоянные. Из равенства  $y(0)=0$ следует $C_1=0$, а из условия $y\in L_2(\bR_+)$ и оценки
\eqref{eq:est5}, которую докажем ниже, получаем $C_2=0$, т.е. $y$ допускает представление \eqref{eq:L-1}. Таким
образом, формула \eqref{eq:L-1} задает обратный оператор $L_c^{-1}$. Его ограниченность следует из оценки
 $|y(x)| \leqslant M\|f\|$ на любом конечном интервале $[0,b]$ и оценки \eqref{eq:est5}.
Здесь и далее через $M$  (или $M_1, M_2$) обозначаются различные положительные константы, а
$\|\cdot\|=\|\cdot\|_{L_2(\bR_+)}$.

Докажем соотношения \eqref{eq:decr}.     Вначале отметим, что в силу \eqref{eq:U1} найдется постоянная $M$, такая, что
при $t>0$  выполняется оценка
$$
|\U(c^{1/3}t)|\le Mt^{-1/4}\exp\left(a t^{3/2}\right), \quad a: = \tfrac 23 |c|^{1/2}\cos\tfrac\gamma2
$$
(напомним, что $\gamma=\arg c$). Легко видеть, что функция $g(t) = t^{-1/4} \exp(at^{3/2})$  возрастает при достаточно
больших $t\ge b= b(a)$. Поэтому при $x> b+1$ имеем
\begin{gather*}
\left|\int_0^x \U(c^{1/3}t)f(t)dt\right|\le \left( \int_0^b + \int_b^{x-1} + \int_{x-1}^x\right)\, |\U(c^{1/3}t)|\,
|f(t)|\,dt\ \le  \\
\le M_1 \|f\|+ M \|f\| x^{1/4} \exp\left[ a(x^{3/2} - x^{1/2})\right] + M x^{-1/4}\exp\left(ax^{3/2}\right)
\left(\int_{x-1}^x|f(t)|^2dt\right)^{1/2}.
\end{gather*}
Здесь $M_1$ --- постоянная, зависящая только от $b$. При переходе ко второму неравенству мы учли, что длина отрезка интегрирования для второго интеграла
меньше  $x$, а также неравенство $(x-1)^{3/2} \le x^{3/2} - x^{1/2}$, верное при достаточно больших $x$.
Из полученной оценки и представления  \eqref{eq:Ai1} при  достаточно больших $x$ имеем
\begin{gather*}
\left|\Ai(c^{1/3}x)\int_0^x\U(c^{1/3}t)f(t)dt\right|\le
Mx^{-1/4}\exp\left(-a x^{3/2}\right)\|f\|+\notag\\
+M\exp\left(-a x^{1/2}\right)\|f\|+Mx^{-1/2}\left(\int_{x-1}^x|f(t)|^2dt\right)^{1/2}.
\end{gather*}
Мы получили  оценку первого слагаемого в \eqref{eq:L-1}. Аналогично получаем оценку второго слагаемого. При достаточно
больших $t$ функция $g(t) = t^{1/2}\exp(-at^{3/2})$  убывает, поэтому при больших $x$
\begin{gather*}
\left|\int_x^\infty
\Ai(c^{1/3}t)f(t)dt\right|\le M\left(\int_{x+1}^\infty g^2(t) t^{-3/2}\, dt\right)^{1/2}\|f\|+ M\int_x^{x+1} g(t)t^{-3/4}\, |f(t|\,dt\\
\le 2^{1/2} M x^{-1/4}g(x+1)\|f\|+ Mg(x)x^{-3/4}\left(\int_x^{x+1}|f(t)|^2dt\right)^{1/2}.
\end{gather*}
Так как $(x+1)^{3/2}\le x^{3/2}+x^{1/2}$, то
$$
\U(c^{1/3}x) \le M x^{1/4} g^{-1}(x), \quad g(x+1) g^{-1}(x) \le M \exp\left(-ax^{1/2}\right)
$$
и модуль второго слагаемого в правой части равенства \eqref{eq:L-1} оценивается величиной
\begin{equation*}
M\exp\left(-ax^{1/2}\right)\|f\| + Mx^{-1/2}\left(\int_{x}^{x+1}|f(t)|^2dt\right)^{1/2}.
\end{equation*}
Cкладывая полученные оценки, приходим к неравенству
\begin{equation}\label{eq:est3}
|y(x)|\le M_2\exp\left(-ax^{1/2}\right)\|f\|+M_2 x^{-1/2}\left(\int_{x-1}^{x+1}|f(t)|^2dt\right)^{1/2},\qquad x>b+1,
\end{equation}
что доказывает первое соотношение в \eqref{eq:decr}. Второе соотношение в \eqref{eq:decr} получается так же, только
вместо равенства \eqref{eq:L-1} используем \eqref{eq:L-1'} и учитываем, что оценки производных $\Ai'$ и $\U'$
отличаются от оценок самих функций множителем $x^{1/2}$.

Докажем, что $x^{1/2}y(x)\in L_2(\bR_+)$. Из оценки \eqref{eq:est3}  получаем
\begin{gather}\label{eq:est5}
\int_b^\infty |x^{1/2}y(x)|^2dx\le M\|f\|^2\int_b^\infty x\exp\left(-2ax^{1/2}\right)dx+
\\+M\int_b^\infty\int_{x-1}^{x+1}|f(t)|^2dt\,dx \le M\|f\|^2+M\int_{b-1}^\infty|f(t)|^2\int_{t-1}^{t+1}dx\,dt\le
M\|f\|^2.\notag
\end{gather}
Включение $y'\in L_2(\bR_+)$ получается аналогично, если учесть, что $|y'(x)|$ оценивается величиной в правой части
\eqref{eq:est3}, умноженной на $x^{1/2}$. Этим завершается доказательство утверждения.
\end{proof}

\begin{Proposition}\label{pr:sect}
Числовой образ оператора $\L_c$ лежит в замкнутом секторе $S_\gamma$ комплексной плоскости, ограниченном лучами $\arg
\la=0$ и $\arg \la=\gamma$, где $\gamma=\arg c$.
\end{Proposition}
\begin{proof}
Квадратичная форма нашего оператора имеет вид
$$
(\L_cy,y)=\int_0^\infty (-y''\overline{y}+cxy\overline{y})dx=-\overline{y}(x)y'(x)\Big|_0^\infty+\int_0^\infty
|y'|^2dx+c\int_0^\infty x|y|^2dx.
$$
Остается заметить, что $y(0)=0$ и $\overline{y}(x)y'(x)\to 0$ при $x\to\infty$ в силу \eqref{eq:decr}.
\end{proof}

\begin{Proposition}\label{pr:conj}
Оператор $\L_c$ замкнут и имеет нулевые дефектные числа (размерность ядра и коразмерность образа). Сопряженный оператор
$\L_c^*$ совпадает с оператором $\L_{\overline c}$.
\end{Proposition}
\begin{proof}
Первое утверждение следует из утверждения \ref{pr:inverse}, поскольку мы предъявили обратный оператор
$\L_c$, который замкнут, поскольку ограничен. Следовательно, обратный к нему также замкнут.  Для
доказательства второго утверждения проверим равенство Лагранжа. Пусть $y(x)\in\Dom(L_c)$, а $u(x)\in\Dom(L_{\overline
c})$. Тогда
\begin{gather*}
(\L_cy,u)=\int_0^\infty\left(-y''(x)\overline u(x)+cy(x)\overline u(x)\right)dx=\int_0^\infty\left(-y(x)\overline
u''(x)+y(x)\overline{\overline c u}(x)\right)dx+\notag\\
+y(x)\overline u'(x)\Big|_0^\infty-y'(x)\overline
u(x)\Big|_0^\infty=\int_0^\infty\left(-y(x)\overline u''(x)+y(c)\overline{\overline c u}(x)\right)dx=(y,\L_{\overline
c}u). 
\end{gather*}
Положим $\L_c y = z, \ \L_{\overline c}u = v$.  Из равенства Лагранжа получаем
$$
(z,\L_{\overline c}^{-1}v) = (\L_c y, u) = (y, \L_{\overline c}u) = (y,v)=(\L_c^{-1} z, v) \quad \forall\ z, v\in
L_2(\mathbb R_+).
$$
Следовательно, $\L_{\overline c}^{-1} = \left(\L_c^{-1}\right)^*$. Но тогда $\L_{\overline c} = \L_c^*$.
\end{proof}

\begin{Proposition}\label{pr:msect}
Резольвента $\R_c(\la)=(\L_c-\la I)^{-1}$ определена и является ограниченным оператором для любого $\la\in\bC\setminus
S_\gamma$, где сектор $S_\gamma$  определен в утверждении \ref{pr:sect}. При этом
\begin{equation}\label{eq:Rest}
\|\R_c(\la)\|_{L_2(\bR_+)}\le\frac1{\dist(\la,S_\gamma)}.
\end{equation}
\end{Proposition}
\begin{proof}
При $\gamma=0$ оператор $\L_c$ самосопряжен и положителен, $S_\gamma=\bR_+$ и доказываемое утверждение хорошо известно.
Пусть $\gamma\ne0$. Согласно утверждению \ref{pr:sect}, оба оператора $T_1=e^{i(\pi/2-\gamma)}\sign\gamma\cdot \L_c$ и
$T_2=-i\sign\gamma\cdot \L_c$ являются аккретивными (т.е. $\Re (T_jy,y)\ge0$ для любого $y\in\Dom(\L_c)=\Dom(T_j)$,
$j=1,\,2$). В силу утверждения
 \ref{pr:conj}, оба этих оператора замкнуты и имеют нулевые дефектные числа, т.е. являются замкнутыми
максимальными аккретивными операторами.
Следовательно
 для любого $z$ из открытой левой полуплоскости операторы $T_1- zI$ и
$T_2-zI$ обратимы и (см., например, \cite[Гл.III.10]{Ka})
$$
\|(T_1-z I)^{-1}\|\le|\Re z|^{-1}\qquad\text{и}\ \|(T_2-z I)^{-1}\|\le|\Re z|^{-1}.
$$
Очевидно, эти оценки эквивалентны оценке \eqref{eq:Rest}.
\end{proof}

Из утверждения \ref{pr:sect} следует, что числовой образ оператора $\L_c$ лежит в секторе $S_\gamma$. Поэтому оператор
$T=e^{-i\gamma/2}\L_c$ является секториальным, а из утверждения \ref{pr:msect} следует, что этот оператор является
$m$-секториальным. С любым $m$-секториальным оператором единственным образом ассоциируется замкнутая секториальная
полуторолинейная форма. Мы найдем ее явный вид.

\begin{Proposition}\label{pr:form}
Замкнутая секториальная полуторалинейная форма $\t$ оператора\break $T=e^{-i\gamma/2}\L_c$ имеет вид
\begin{gather}
\t[u,v]=e^{-i\gamma/2}\int_0^\infty u'(x)\overline v'(x)\,dx+|c|e^{i\gamma/2}\int_0^\infty x u(x)\overline v(x)\,dx,\notag \\
u,\,v\in\Dom(\t)=\left\{y\in L_2(\bR_+):y',\,x^{1/2}y\in L_2(\bR_+)\right\}. \label{eq:tform2}
\end{gather}
\end{Proposition}

\begin{proof}
Определим форму $\t_0[u,v]=(Tu,v)$ на области $\Dom(T)=\Dom(\L_c)$. Интегрируя по частям, так же как в доказательстве
утверждения \ref{pr:sect}, получаем
$$
\t_0[u,v]=e^{-i\gamma/2}\int_0^\infty u'(x)\overline v'(x)\,dx+|c|e^{i\gamma/2}\int_0^\infty xu(x)\overline v(x)\,dx.
$$
Очевидно, определенный выше линеал $\Dom(\t)$ является замкнутым множеством в норме, порожденной скалярным
произведением
$$
[u,v]=\int_0^\infty\, (xu(x)\overline v(x)+u'(x)\overline v'(x))\,dx,
$$
т.е. $\Dom(\t)$  является гильбертовым пространством с таким скалярным произведением.  Заметим, что $\Dom(\L_c)$
является ядром формы $\t$ (см. \cite[Гл. VI, Теорема 2.1]{Ka}). Поэтому область определения  замыкания формы $\t_0$
совпадает с $\Dom(\t)$, и форма $\t$, отвечающая оператору  $T$, имеет вид \eqref{eq:tform2}.
\end{proof}

\begin{Proposition}\label{pr:comp}
Для любого $\la$ из резольвентного множества оператор $\R_c(\la)$ компактен. Спектр оператора $\L_c$ дискретен и
состоит из последовательности простых собственных значений
\begin{equation}\label{eq:spectr}
\la_n=t_nc^{2/3},\ n\in\bN,\quad \text{где}\quad t_n>0\quad \text{и}\quad
t_n=\left[(3\pi/2)\left(n-1/4\right)\right]^{2/3}\!\!\!+O(n^{-4/3}).
\end{equation}
 Им отвечают собственные функции
 \begin{equation}\label{eq:y_n}
  y_n(x) = \Ai(-t_n + xc^{1/3}), \quad n = 1, 2. \dots .
 \end{equation}
\end{Proposition}

\begin{proof}
Компактность резольвенты оператора $T=e^{-i\gamma/2}\L_c$ эквивалентна компактности резольвенты оператора
 $H=\Re T$  (см. \cite[Гл. VI, Теорема 3.3]{Ka}).  Оператор $H$ порождается  квадратичной формой
 $\h=\Re\t$, т.е.
$$
\h[u,v]=\cos(\tfrac\gamma2)\int_0^\infty u'(x)\overline v'(x)\,dx+|c|\cos(\tfrac\gamma2)\int_0^\infty xu(x)\overline
v(x)\,dx,\qquad\Dom(\h)=\Dom(\t).
$$
Оператор $H$ однозначно восстанавливается по квадратичной форме, поэтому из утверждения
\ref{pr:form} получаем
\begin{gather}
Hy=\cos(\tfrac\gamma2)(-y''+|c|xy),\label{eq:H}\\
\Dom(H)=\left\{y\in L_2(\bR_+):y\in W_{2,loc}^2,\ -y''+|c|xy\in L_2(\bR_+),\ y(0)=0\right\}.\notag
\end{gather}
Компактность резольвенты  $(H-\lambda)^{-1}$  следует из критерия Молчанова (см.,
например, \cite[Гл. VII, \S 24]{Na}). Итак, резольвента $\R_c(\la)$ компактна, а спектр $\sigma(\L_c)$ дискретен.
Запишем уравнение на собственные значения
$$
-y''+cxy=\la y\quad y(0)=0,\ y\in L_2(\bR_+).
$$
Сделаем замену $t=-\la c^{-2/3}+xc^{1/3}$. Тогда уравнение можно записать в виде $-y''_{tt}+ty=0$, т.е.
\begin{equation}\label{eq:EV1}
y(x)=C_1\Ai(c^{1/3}x-\la c^{-2/3})+C_2\U(c^{1/3}x-\la c^{-2/3}).
\end{equation}
Учитывая \eqref{eq:Ai1} и \eqref{eq:U1}, получаем $C_2=0$, а поскольку $y(0)=0$, приходим к уравнению на собственные
значения $\Ai(-\la c^{-2/3})=0$. Остается заметить, что все нули функции $\Ai(z)$ просты, а равенства \eqref{eq:spectr}
теперь следуют из \eqref{eq:Aizero}. Подставляя в \eqref{eq:EV1} $\lambda = \lambda_n$, получаем равенства
\eqref{eq:y_n}.
\end{proof}

\section*{Теорема о полноте для оператора $\L_c$}

\begin{Theorem}\label{tm:56}
Система собственных функций оператора $\L_c$ полна и минимальна в $L_2(\bR_+)$ при условии $|\arg c|<(5\pi)/{6}$.
\end{Theorem}

\begin{proof}
Минимальность системы собственных функций  $\{y_n\}$  оператора $\L_c$ следует из известных соотношений
$$
(y_n, z_k) = c_n \, \delta_{nk},\ \quad c_n \ne 0,
$$
где $\{z_k\}$ --- система собственных функций сопряженного оператора $\L_{\overline c}$. Для доказательства полноты
будем использовать метод Левинсона  (см. \cite[Приложение 4]{Le}), принимая во внимание наличие функции
$\Ai(\w+c^{1/3}x)$, порождающей собственные функции  $y_n$ при $\w=-t_n$. Для определенности далее будем рассматривать
случай $\Im c\ge0$, т.е. $\gamma\in[0,\pi)$, где $\gamma=\arg c$. Случай $\gamma\in(-\pi,0]$ рассматривается аналогично
(проводимые ниже оценки в секторах комплексной плоскости получаются аналогично в секторах, симметричных относительно
вещественной оси).

Доказательство разобьем на несколько этапов. Пусть функция $f\in L_2(\mathbb R_+)$ ортогональна собственным функциям
оператора  $\L_c$. Рассмотрим функцию
\begin{equation}\label{eq:Fdef}
F(\w)=\frac{F_0(\w)}{\Ai(\w)},\qquad F_0(\w)=\int_0^\infty \Ai(\w+xc^{1/3})f(x)dx.
\end{equation}
На первом шаге мы покажем, что функция $F$
 является целой функцией порядка $\rho \le 3/2$ и конечного типа при $\rho = 3/2$. На втором шаге покажем, что
 функция $F$ допускает оценку
 \begin{equation}\label{eq:F}
 |F(\w)| \le M\|f\| R^{1/2}, \qquad R=|\w|\ge 1
 \end{equation}
в секторе
 \begin{equation}\label{eq:S}
 S= \left\{\w\in \mathbb C: \ -\pi +\gamma/{3} \le \arg\w \le \pi -{2\gamma}/3 \right\}.
 \end{equation}
На третьем шаге мы покажем,
что существует число $\alpha_0 \in (0, \gamma/3)$, такое, что оценка \eqref{eq:F} остается справедливой
в секторе $S' = S\bigcup S_0$, где
\begin{equation}\label{eq:S_0}
 S_0= \left\{\w\in \mathbb C: \ -\pi  +\alpha_0  \le \arg \w \le -\pi +\gamma/{3} \right\}.
 \end{equation}
Число $\alpha_0$  не удается вычислить явно (оно является корнем трансцендентного уравнения), но удается показать, что
раствор дополнительного сектора  $\mathbb C \setminus S'$ меньше  $2\pi/3$, если только $\gamma < 5\pi/6$.  Это влечет
выполнение оценки \eqref{eq:F} во всей комплексной плоскости. На четвертом шаге  мы покажем  $F(\w) \equiv 0$, а на
пятом шаге получим $f(x) \equiv 0$. Приступим к реализации намеченного плана.

\medskip\noindent\textit{Шаг 1.} Покажем, что функция $F_0$ в \eqref{eq:Fdef}   корректно определена и
голоморфна по параметру $\w\in\bC$. Зафиксируем
произвольное число $\delta\in(0,\pi/4]$, такое, что $\gamma/2+\delta<\pi/2$. Пусть $\w$ пробегает компакт $|\w|\le R$.
Представим  функцию $F_0$  в виде
$$
F_0(\w)=\left(\int_0^{x_0}+\int_{x_0}^\infty\right) \Ai(\w+xc^{1/3})f(x)dx:=F_1(\w)+F_2(\w),\ \ x_0 =
x_0(R)=\frac{R}{|c|^{1/3}\sin(2\delta/3)}.
$$
Первый интеграл является собственным и, стало быть, голоморфен по параметру $\w$ в круге $|\w|\le R$, а для
доказательства голоморфности второго интеграла заметим, что при $x> x_0(R)$
\begin{equation*}
\left|\arg\left(1+\frac{\w}{c^{1/3}x}\right)\right|\le\frac{2\delta}3,\quad   \Longrightarrow\ \
\arg(\w+c^{1/3}x)^{3/2}\in\left[\frac\gamma2-\delta,\frac\gamma2+\delta\right].
\end{equation*}
Тогда при $x> x_0(R)$ и $|\w|\le R$ имеем
\begin{equation}\label{eq:a_1}
\Re(\w +c^{1/3}x)^{3/2} \ge \left|\w +c^{1/3}x\right|^{3/2} \cos \left(\frac\gamma{2}+\delta\right) \ge \frac 32\, a_1
\left( |c|^{1/3}x -R\right)^{3/2},  .
\end{equation}
где $a_1: =2/3\cos\left(\gamma/{2}+\delta\right)$. Следовательно, в силу \eqref{eq:Ai1} и \eqref{eq:a_1} при $x>x_0(R)$
\begin{equation}\label{eq:F2est}
|\Ai(\w+c^{1/3}x)|\le M\exp\left\{\tfrac23\Re(\w+c^{1/3}x)^{3/2}\right\}\le
M\exp\left\{-a_1\left(|c|^{1/3}x-R\right)^{3/2}\right\}.
\end{equation}
Аналогичная оценка сохраняется и для $\Ai'(\w+c^{1/3} x)^{3/2}$, только правую часть нужно умножить на
$|\w+c^{1/3}x|^{1/4}$. Отсюда следует, что интеграл $F_2(\w)$  сходится равномерно по параметру $\w$  в круге  $|\w|\le
R$ и остается равномерно сходящимся после дифференцирования по $\w$, т.е.   функция $F_2$  голоморфна в круге $|\w|\le
R$. Поскольку $R$ произвольно, функция $F_0 =F_1 +F_2$  является целой. Более того, целой является и функция $F(\w)$,
поскольку все точки $\w=-t_k$ являются ее устранимыми особенностями (по нашему предположению, функция $F_0(\w)$
обращается в этих точках в нуль, а функция $\Ai(\w)$ имеет в этих точках простые нули).

Теперь оценим рост функции $F_0(\w)$ при $|\w|\to\infty$. Длина отрезка интегрирования  в интеграле $F_1$  пропорциональна  $R$, поэтому
\begin{gather}\label{eq:F1est}
|F_1(\w)|\le MR^{1/2}\|f\|\exp\left\{\tfrac23(R+|c|^{1/3}x_0(R))^{3/2}\right\}\le M R^{1/2} \|f\| \exp\left(M_1 R^{3/2}\right) .
\end{gather}
Положим $v(x) = \exp \left\{-a_1(|c|^{1/3}x - R)^{3/2} +x\right\}$, где число  $a_1$  определено в \eqref{eq:a_1}.
Эта функция убывает при $x> x_0$, если число $x_0 =x_0(R)$ достаточно большое. Поэтому из \eqref{eq:F2est} получаем
\begin{gather}\label{eq:F2est1}
 |F_2(\w)|\le
M\int_{x_0}^\infty \exp\left\{-a_1(|c|^{1/3}x-R)^{3/2}+x\right\} \, |f(x)| e^{-x}\, dx\\
\le 2^{-1/2}M \|f\|\exp\left\{-a_1(|c|^{1/3}x_0(R)-R)^{3/2}+x_0)\right\}, \quad x_0 =\frac  R{|c|^{1/3}\sin (2\delta/3)}. \notag
\end{gather}
Следовательно, $|F_2(\w)| \to 0$  при $\w \to \infty$.  Из \eqref{eq:F1est} следует, что $F_0 = F_1 +F_2$ --- целая
функция порядка не выше $3/2$ и конечного типа при порядке $3/2$. Тогда целая функция $F(\w)$, равная отношению двух
целых функций $F_0(\w)$ и $\Ai(\w)$, имеет характеристику роста, не превосходящую максимальную характеристику роста
числителя и знаменателя (см., например, \cite[Гл. I, \S 9]{Le}), т.е. функция $F$ имеет порядок $\rho\le3/2$, а при
$\rho=3/2$, конечный тип.

\medskip\noindent\textit{Шаг 2.} Рассмотрим секторы
\begin{equation}\label{eq:S_1}
\begin{split}
S_1 =\left\{\w\in\mathbb C:\ -\pi+\gamma/3 \le \arg\w \le -\pi/3 +\varepsilon\right\},\\ S_2 =\left\{\w\in\mathbb C:\
\pi/3 -\varepsilon \le\arg\w \le \pi- {2\gamma}/3\right\},
\end{split}
\end{equation}
где число $\varepsilon\in(0,\pi/3)$ выберем ниже. Согласно \eqref{eq:F2est1} и \eqref{eq:Ai1},
\begin{equation}\label{eq:F2Ai}
\frac{|F_2(\w)|}{|\Ai(\w)|}\le
M\|f\|\exp\left\{\frac23\left(\cos(\tfrac{3\varphi}2)-\frac{\cos(\tfrac\gamma2+\delta)(1-\sin(\tfrac{2\delta}3))^{3/2}}{\sin^{3/2}(\tfrac{2\delta}3)}\right)R^{3/2}+
\frac{R}{|c|^{1/3}\sin(\tfrac{2\delta}3)}\right\},
\end{equation}
где $\varphi=\arg\w$. В объединении $S_1\bigcup S_2$ секторов, определенных в  \eqref{eq:S_1}, выполнена оценка
$\cos(3\varphi/2)\le\sin(3\eps/2)$. Выберем число  $\eps>0$ так, чтобы
$$
\sin\left(\frac{3\eps}2\right)-\frac{\cos(\tfrac\gamma2+
\delta)(1-\sin(\tfrac{2\delta}3))^{3/2}}{\sin^{3/2}(\tfrac{2\delta}3)}<0.
$$
Тогда при  $R\to\infty$, получаем оценку
\begin{equation}\label{eq:o}
|F_2(\w)|/|\Ai(\w)|=o(1)\|f\|.
\end{equation}
Далее, вновь применяя \eqref{eq:Ai1}, получаем
\begin{equation}\label{eq:F1Ai1}
|F_1(\w)|/{|\Ai(\w)|}\le MR^{1/2}\|f\|\max_{x>0}e^{\xi(x)},
\end{equation}
где
\begin{equation}\label{eq:xidef}
\xi(x)=\frac23\left(\cos(\tfrac{3\varphi}2)R^{3/2}-\Re(\w+c^{1/3}x)^{3/2}\right), \quad x\in[0,+\infty).
\end{equation}
Заметим, что $\xi'(x)=-\Re\left(c^{1/3}(\w+c^{1/3}x)^{1/2}\right)$, откуда следует, что $\xi'(x)\le0$, если
\begin{equation}\label{eq:xi'}
\arg(\w+c^{1/3}x)\in[-\pi,\pi-2\gamma/3].
\end{equation}
В этом случае $\xi(x)$ монотонно невозрастает, а так как $\xi(0)=0$, то $\xi(x)\le0$. Поэтому при выполнении условия
\eqref{eq:xi'} из \eqref{eq:F1Ai1} получаем
\begin{equation}\label{eq:F1Ai}
|F_1(\w)|/{|\Ai(\w)|}\le M\|f\|R^{1/2}\quad\text{при}\ \  |\w|=R\to\infty.
\end{equation}
Проверим,  что условие \eqref{eq:xi'} выполнено в объединении $S_1\bigcup S_2$   для любого $x\ge0$. Разберем два
случая. Если $\varphi\in[-\pi+\gamma/3,-\pi/3+\eps]$, то оба вектора $\w$ и $c^{1/3}x$ лежат в полуплоскости
$-\pi+\gamma/3\le\arg z\le\gamma/3$, а значит в той же полуплоскости лежит их сумма. Если же
$\varphi\in[\pi/3-\eps,\pi-2\gamma/3]$, то оба вектора $\w$ и $c^{1/3}x$ лежат в секторе $0\le\arg z\le \pi-2\gamma/3$
раствора $\le\pi$, а значит в том же секторе лежит их сумма. В обоих случаях условие \eqref{eq:xi'} выполнено и мы
доказали, что в объединении $S_1\bigcup S_2$ секторов, заданных \eqref{eq:S_1}, справедлива оценка \eqref{eq:F1Ai}.
Учитывая оценку \eqref{eq:o}, получаем, что в объединении $S_1\bigcup S_2$ справедлива асимптотическая оценка
\begin{equation}\label{eq:Finest}
|F(\w)|\le M\|f\|R^{1/2}.
\end{equation}
Сектор $S$ определен  формулой \eqref{eq:S}, а  раствор сектора $S\setminus(S_1\bigcup S_2)$ меньше $2\pi/3$. Применив
теорему Фрагмена--Линделефа к функции  $\widetilde F(\w)=F(\w)(\w+1)^{-1/2}$, получим, что оценка \eqref{eq:Finest}
выполняется во всем секторе $S$.

\medskip\noindent\textit{Шаг 3.} Остается оценить функцию $F$  в секторе $S'' =\mathbb C\setminus S$. Раствор
сектора $S''$ равен $\gamma$; поэтому оценка $F$  в оставшемся секторе также получается из теоремы Фрагмена--Линделефа,
если $\gamma <2\pi/3$. В случае $\gamma\in[2\pi/3,5\pi/6)$, который мы рассмотрим теперь, требуются дополнительные
оценки.

Пусть $-\pi<\varphi<-\pi+\gamma/3$. Оценка дроби $|F_2(\w)|/|\Ai(\w)|$ не меняется (эта дробь по-прежнему ограничена
некоторой константой), поскольку $\cos(3\varphi/2)<0$ и мы можем воспользоваться оценкой \eqref{eq:F2Ai}. Докажем
оценку \eqref{eq:F1Ai} в секторе $S_0$, определенном \eqref{eq:S_0}, где $\alpha_0 <\gamma/3$. Исследуем на максимум
функцию $\xi(x)$, определенную в \eqref{eq:xidef}. Напомним, что ветвь отображения $z^{3/2}$ зафиксирована нами выбором
аргумента $\arg z\in[-\pi,\pi)$. Легко видеть, что при $\varphi\in(-\pi,-\pi+\gamma/3)$ луч $\w+c^{1/3}x$, $x>0$,
пересекает луч $(-\infty,0)$, а значит кривая $(\w+c^{1/3}x)^{3/2}$, $x\in[0,+\infty)$, имеет разрыв первого рода в
точке, которую мы обозначим $x_1$. Заметим, однако, что величина скачка в точке $x_1$ --- чисто мнимое число, так что
функция $\xi(x)$ непрерывна. Проведем исследование функции $\xi$ отдельно на промежутках $x\in(0,x_1)$ и
$x\in(x_1,+\infty)$. На первом промежутке $\xi'(x)<0$, так как аргумент числа $\w+c^{1/3}x$ меняется от $\varphi$ (при
$x=0$) до $-\pi$ (при $x=x_1$), т.е. выполняется условие \eqref{eq:xi'}. На втором промежутке аргумент числа
$\w+c^{1/3}x$ монотонно убывает от $\pi$ (при $x=x_1+0$) до $\gamma/3$ (при $x\to+\infty$), а значит найдется
единственная точка, которую мы обозначим $x_2$, в которой $\arg(\w+c^{1/3}x)=\pi-2\gamma/3$. В точке $x_2$ производная
$\xi'$ меняет знак с плюса на минус, т.е. $x_2$
--- точка локального максимума функции $\xi$. Таким образом, оценка \eqref{eq:F1Ai}, которая, в силу \eqref{eq:F1Ai1} эквивалентна
условию $\xi(x)\le0$ при $x\in[0,+\infty)$, выполняется в точности тогда, когда $\xi(x_2)\le0$. Чтобы вычислить
значение $\xi(x_2)$ рассмотрим на комплексной плоскости треугольник с вершинами  в точках $0$, $\w$ и $\w+c^{1/3}x_2$.
Соответствующие углы этого треугольника равны
 $\al+2\gamma/3$, $\gamma/3-\al$ и $\pi-\gamma$ (здесь мы ввели
обозначение $\al=\pi+\varphi\in(0,\gamma/3)$),\ $\varphi =\arg\w$. В силу теоремы синусов
\begin{multline*}
|\w+c^{1/3}x_2|=R\frac{\sin(\gamma/3-\al)}{\sin\gamma} \\ \Longrightarrow\
\Re(\w+c^{1/3}x_2)^{3/2}=|\w+c^{1/3}x_2|^{3/2}\cos(\tfrac{3\pi}2-\gamma)=
-R^{3/2}\frac{\sin^{3/2}(\tfrac\gamma3-\al)}{\sin^{1/2}\gamma},\\
\xi(x_2)=\frac23R^{3/2}\left(\cos(\tfrac{3\varphi}2)+\frac{\sin^{3/2}(\tfrac\gamma3-\al)}{\sin^{1/2}\gamma}\right)=
\frac23R^{3/2}\left(\frac{\sin^{3/2}(\tfrac\gamma3-\al)}{\sin^{1/2}\gamma}-\sin(\tfrac{3\al}2)\right).
\end{multline*}
Итак, неравенство $\xi(x_2)\le0$ равносильно условию
\begin{equation}\label{eq:ineqal}
\sin^{3/2}(\gamma/3-\al)\sin^{-1/2}\gamma-\sin({3\al}/2)\le0.
\end{equation}
Обозначим правую часть этого неравенства --- функцию переменной $\al\in[0,\gamma/3]$ с параметром
$\gamma\in[2\pi/3,\pi)$ через $\eta(\al)$. Легко видеть, что $\eta(\al)$ монотонно убывает, $\eta(0)>0$ и
$\eta(\gamma/3)<0$. Таким образом, неравенство \eqref{eq:ineqal} выполнено на отрезке $\al\in[\al_0,\gamma/3]$ для
некоторого $\al_0\in(0,\gamma/3)$. Это означает, что оценка \eqref{eq:F1Ai}, а значит и \eqref{eq:Finest}, доказана
нами для всех лучей $\arg\w=\varphi\in[-\pi+\al_0,\pi-2\gamma/3]$. Принцип Фрагмена--Линделефа позволит нам
распространить эту оценку на всю комплексную плоскость, если раствор оставшегося угла, равный $\al_0+2\gamma/3$ строго
меньше $2\pi/3$. Число $\al_0=\al_0(\gamma)$ является корнем трансцендентного уравнения $\eta(\al)=0$. Мы не будем его
искать в каком бы то ни было виде, а вместо этого заметим, что в силу монотонности функции $\eta$

$$
\al_0+\frac{2\gamma}3<\frac{2\pi}3\ \ \Longleftrightarrow\ \ \al_0<\frac{2(\pi-\gamma)}3\ \ \Longleftrightarrow\ \
\eta\left(\frac{2(\pi-\gamma)}3\right)<0
$$

$$
\Longleftrightarrow\ \ \sin^{3/2}
(\gamma-\tfrac{2\pi}3) \sin^{-1/2}\gamma-
\sin\gamma <0\ \ \Longleftrightarrow\ \
\sin^{3/2}(\gamma-\tfrac{2\pi}3)<\sin^{3/2}\gamma\ \ \Longleftrightarrow\ \ \gamma<\frac{5\pi}6.
$$

\medskip\noindent\textit{Шаг 4.} Итак, при $\gamma<5\pi/6$ целая функция $F(\w)$ допускает асимптотическую оценку
\eqref{eq:Finest} во всей комплексной плоскости, а значит, равна константе. Докажем, что эта константа равна нулю. Для
этого достаточно убедиться, что $F(\w)\to0$ при $\w\to\infty$ хотя бы по одному лучу в комплексной плоскости. Вернемся
к началу шага 2, выберем луч $\w=Re^{-i\pi/2}$ и заметим, что он попадает в сектор $S_1$, а потому выполнена оценка
\eqref{eq:o}. Остается усилить оценку \eqref{eq:F1Ai}. Для этого разобьем путь интегрирования $x\in[0,x_0]$ на два
участка точкой $x_3(R)=R^{-\theta}$, где $\theta\in(0,1/2)$ произвольно. Так как
$$
\frac{|F_1(\w)|}{|\Ai(\w)|}\le M\int_0^{x_0(R)}|f(x)|e^{\xi(x)}\,dx,
$$
где функция $\xi(x)$ определена в \eqref{eq:xidef} и монотонно убывает на $[0,+\infty)$ от $\xi(0)=0$ к $-\infty$, то
$$
\left(\int_0^{x_3}+\int_{x_3}^{x_0}\right)|f(x)|e^{\xi(x)}\,dx\le
\|f\|\left(\sqrt{x_3}+\sqrt{x_0}e^{\xi(x_3)}\right)\le \|f\|\left(R^{-\theta/2}+CR^{1/2}e^{\xi(x_3)}\right).
$$
Остается заметить, что при $R\to\infty$
$$
\xi(x_3(R))\sim-\Re(c^{1/3}x_3\w^{1/2})=-|c|^{1/3}R^{-\theta+1/2}\cos\left(\gamma/3-\pi/4\right),
$$
а значит $R^{1/2}e^{\xi(x_3(R))}=o(1)$. Следовательно, $F(\w) \equiv 0.$

\medskip\noindent\textit{Шаг 5.}
Мы доказали, что функция $F_0(\w)=\int_0^\infty \Ai(\w+c^{1/3}x)\,d\mu \equiv 0$ (здесь и далее для сокращения записи
обозначаем $d\mu=f(x)dx$). Напомним, что $\Ai''(t)=t\Ai(t)$, откуда по индукции
\begin{equation*}
\Ai^{(n)}(t)=P_n(t)\Ai(t)+Q_n(t)\Ai'(t),
\end{equation*}
где $P_n$ и $Q_n$ --- многочлены. При этом
\begin{align*}
P_0(t)&=1,\quad & P_1(t)&=0,\quad & P_2(t)&=t,\quad & P_3(t)&=1,\quad & P_4(t)&=t^2,\quad & P_5(t)&=4t,\\
Q_0(t)&=0 & Q_1(t)&=1, & Q_2(t)&=0, & Q_3(t)&=t, & Q_4(t)&=2, & Q_5(t)&=t^2,
\end{align*}
$$
P_n(t)=P'_{n-1}(t)+tQ_{n-1}(t), \qquad Q_n(t)=P_{n-1}(t)+Q'_{n-1}(t).
$$
Положим $\deg P_n(t)=p_n$, $\deg Q_n(t)=q_n$. Тогда для всех $n\ge5$
$$
\begin{cases}p_n=q_{n-1}+1,\\q_n=p_{n-1}\end{cases}\Longrightarrow\begin{cases}p_n=p_{n-2}+1,\\q_n=q_{n-2}+1\end{cases},
$$
откуда легко выводим, что для всех $n\ge 2$
\begin{align}\notag
p_{2n}&=n,\quad & p_{2n+1}&=n-1,\\ q_{2n}&=n-2, & q_{2n+1}&=n.\label{eq:degs}
\end{align}
Проводя дифференцирование по переменной $\w$, получим
\begin{equation}\label{eq:gj}
F_0^{(j)}(\w)=\!\!\int_0^\infty\!\!\!\! P_j(\w+c^{1/3}x)\Ai(\w+c^{1/3}x)+Q_j(\w+c^{1/3}x)\Ai'(\w+c^{1/3}x)\,d\mu\equiv0
\end{equation}
для всех $j\ge0$. В частности,
\begin{gather*}
F_0(\w)=\int_0^\infty \Ai(\w+c^{1/3}x)\,d\mu\equiv0,\qquad F_0'(\w)=\int_0^\infty \Ai'(\w+c^{1/3}x)\,d\mu\equiv0,\\
F_0''(\w)=\int_0^\infty(\w+c^{1/3}x)\Ai(\w+c^{1/3}x)\,d\mu\equiv0\ \ \Longrightarrow\ \  \int_0^\infty x\Ai(\w+c^{1/3}x)\,d\mu\equiv0,\\
\begin{split}
F_0'''(\w)=\int_0^\infty (\Ai(\w+c^{1/3}x)+(\w+c^{1/3}x)\Ai'(\w+c^{1/3}x))\,d\mu\equiv0\\ \Longrightarrow\ \
\int_0^\infty x\Ai'(\w+c^{1/3}x)\,d\mu\equiv0.
\end{split}
\end{gather*}
Докажем по индукции, что
\begin{equation}\label{eq:powers}
\int_0^\infty x^n\Ai(\w+c^{1/3}x)\,d\mu\equiv\int_0^\infty x^n\Ai'(\w+c^{1/3}x)\,d\mu\equiv0,\quad n=0,\,1,\,2,\dots.
\end{equation}
База индукции --- случай $n=0$ и $n=1$. В этом случае равенства уже доказаны. Если равенства \eqref{eq:powers} доказаны
для всех $n< N$, где $N\ge2$, то, записав \eqref{eq:gj} для $j=2N$ и учитывая \eqref{eq:degs}, получим
\begin{multline*}
\int_0^\infty\left((\w+c^{1/3}x)^N+\sum_{k=0}^{N-1}a_{j,k}(\w+c^{1/3}x)^k\right)\Ai(\w+c^{1/3}x)+\\
+\left(\sum_{k=0}^{N-2}b_{j,k}(\w+c^{1/3}x)^k\right)\Ai'(\w+c^{1/3}x)\,d\mu\equiv0,
\end{multline*}
откуда, в силу предположения индукции,
$$
\int_0^\infty x^N\Ai(\w+c^{1/3}x)\,d\mu\equiv0.
$$
Записав теперь равенство \eqref{eq:gj} для $j=2N+1$ и вновь приняв во внимание \eqref{eq:degs}, получим
\begin{multline*}
\int_0^\infty\left(\sum_{k=0}^{N-1}a_{j,k}(\w+c^{1/3}x)^k\right)\Ai(\w+c^{1/3}x)+\\
+\left((\w+c^{1/3}x)^N+\sum_{k=0}^{N-1}b_{j,k}(\w+c^{1/3}x)^k\right)\Ai'(\w+c^{1/3}x)\,d\mu\equiv0,
\end{multline*}
откуда, в силу предположения индукции и уже проведенного шага,
$$
\int_0^\infty x^N\Ai'(\w+c^{1/3}x)\,d\mu\equiv0.
$$
Равенства \eqref{eq:powers} доказаны.

Теперь из равенств \eqref{eq:powers} стандартным приемом получим $f(x) \equiv 0$. Рассмотрим преобразование Фурье
$$
h(\la)=\int_0^\infty e^{-i\la x}\Ai(c^{1/3}x)f(x)\,dx.
$$
В силу \eqref{eq:Ai1} эта функция является целой функцией переменного $\la$. Заметим, что
$$
h^{(n)}(0)=(-i)^n\int_0^\infty x^n\Ai(c^{1/3}x)f(x)\,dx=0,\ \ n=0,\,1,\,\dots,
$$
а значит $h(\la)\equiv0$. Тогда в силу инъективности преобразования Фурье тождественно нулевой является и функция
$\Ai(c^{1/3}x)f(x)$, поэтому  $f(x)\equiv 0$. Этим завершается доказательство теоремы 1.
\end{proof}

\section*{Теорема о полноте для оператора $L$}

\begin{Theorem}
Пусть оператор $L$ задан дифференциальным выражением \eqref{eq:main2} на области
$$
\Dom (L) = \{y\in L_2(\mathbb R_+): \, y, y' \in AC_{loc} \ \, l(y) \in L_2(\mathbb R_+), \ y(0) =0\},
$$
где  $AC_{loc}$ ---  пространство локально абсолютно непрерывных функций, и пусть выполнены условия \eqref{eq:cond}.
Если
$$
\arg (c_0+i) =:\gamma < 2\pi\alpha/(2+\alpha),
$$
то система корневых функций оператора L полна в $L_2(\mathbb R_+)$. Более того, эта система образует базис для метода
суммирования Абеля порядка $\beta$  при любом  $\beta \in \left( \frac{2+\alpha}{2\alpha}, \frac\pi\gamma, \right).$
\end{Theorem}

\begin{proof} Сначала поясним, что означает понятие базиса для суммирования методом Абеля. Для простоты считаем, что все
собственные значения простые (в общем случае см. \cite[$\S 5$]{Shka1}). Рассмотрим
ряд
$$
S(t, f) = \sum_{k=0}^\infty\, \exp\left\{ \left(-e^{-i\gamma/2}\lambda_k\right)^\beta t \right\}\,
(f, z_k)y_k,
$$
где $\{y_k\}$ ---   система собственных функций оператора $L$, отвечающая собственным значениям $\lambda_k$, а
$\{z_k\}$ --- биортогональная к ней система. Второе утверждение  теоремы означает, что $\forall f \in L_2(\mathbb R_+)$
ряд $S(t, f) $  сходится  $\forall t>0$ по норме пространства $L_2(\bR_+)$ и существует сильный предел  $S(t,f) \to f$
при $t\to +0$.  Конечно, если система $\{y_k\}$ есть базис для метода суммирования Абеля, то она полна, так как из
определения следует, что любая функция $f$ может быть приближена конечными линейными комбинациями системы с
произвольной точностью.

Не ограничивая общности считаем, что постоянные $M_0,\, M_1$ в условиях \eqref{eq:cond} положительны. Пусть $\Dom_0(L)$
подмножество в  $\Dom(L)$, состоящее из функций с компактным носителем. Повторяя рассуждения из работы Лидского
\cite{Li1},  получаем, что в условиях теоремы имеет место случай точки Вейля, т.е. только одно решение уравнения  $l(y)
=0$ принадлежит $L_2(\bR_+)$. В этом случае оператор $L$ имеет ограниченный обратный. Более того, в этом случае
$\Dom_0(L)$ является ядром оператора $L$, т.е. замыкание сужения оператора $L$ на  $\Dom_0(L)$ совпадает с $L$. Но для
всех $f\in \Dom_0(L)$ числовой образ  $(Lf,f)$  лежит в секторе, ограниченным в верхней полуплоскости лучами
$\arg\gamma$  и $\mathbb R_+$. Это утверждение вытекает из  \eqref{eq:cond}   после интегрирования по частям. Из
обратимости $L$ следует, что оператор $T=e^{-i\gamma/2}L$  является $m$-секториальным с раствором  угла $\gamma$. Тогда
оператор $H=\Re T$ задается на области  $\Dom(H) =\Dom(L)$  дифференциальным выражением
$$
l(y) = -y'' +\widetilde r(x) y,
$$
причем из условий \eqref{eq:cond} следует оценка  $\widetilde r(x) \ge Mr(x) \ge M a x^\alpha$.
В работе Лидского доказано, что  при выполнении такого неравенства собственные значения $s_n$  оператора $H$ подчинены оценке
$$
s_n \ge M\, n^{2\alpha/(2+\alpha)}, \qquad n=1, 2, 3,\dots ,
$$
т.е. порядки операторов $H$ и $L$ не меньше  $2\alpha/(2+\alpha)$.  Теперь утверждение теоремы следуют из теоремы
Лидского--Мацаева, см. \cite[Теорема 5.1]{Shka1}.
\end{proof}
\begin{Note} В случае четного потенциала полнота собственных функций оператора Шредингера
на полуоси с условиями Дирихле и Неймана соответственно, влечет за собой полноту собственных функций этого оператора на
всей оси. Действительно, нечетные и четные продолжения собственных функций операторов Дирихле и Неймана  на полуоси
задают собственные функции оператора на всей оси.
\end{Note}

\end{document}